\documentclass[11pt]{article}

\usepackage{amsfonts}
\usepackage{mathrsfs}
\usepackage{bbm}
\usepackage{CJK}
\usepackage{amssymb,amsmath}
\usepackage[dvips]{color}
\usepackage[active]{srcltx}

\allowdisplaybreaks

\newtheorem{lemma}{Lemma}[section]
\newtheorem{theorem}{Theorem}[section]
\newtheorem{corollary}{Corollary}[section]

\newtheorem{proposition}{Proposition}[section]

\def\blemma{\begin{lemma}\it{}\def\elemma{\end{lemma}}}
\def\btheorem{\begin{theorem}\it{}\def\etheorem{\end{theorem}}}
\def\bcorollary{\begin{corollary}\it{}\def\ecorollary{\end{corollary}}}
\def\bproposition{\begin{proposition}\it{}\def\eproposition{\end{proposition}}}

\def\beqlb{\begin{eqnarray}}\def\eeqlb{\end{eqnarray}}
\def\beqnn{\begin{eqnarray*}}\def\eeqnn{\end{eqnarray*}}

\def\qed{\hfill$\Box$\medskip}

\begin{document}

\bigskip

\bigskip

\bigskip

\title{{\bf Branching structure for the transient $(1,
R)$--random walk in random environment  and its applications  }
\thanks{\,The project is partially supported by
National Nature Science Foundation of China (Grant No.\,10721091)
and NCET (No.\,05-0143).}}

\author{Wenming Hong\footnote{School of Mathematical Sciences
\& Laboratory of Mathematics and Complex Systems, Beijing Normal
University, Beijing 100875, P.R.China.
 Email: wmhong@bnu.edu.cn} \ \ \ and  \ \ \ Lin Zhang \footnote{School of Mathematical Sciences
\& Laboratory of Mathematics and Complex Systems, Beijing Normal
University, Beijing 100875, P.R.China.
 Email:zhanglin2008@mail.bnu.edu.cn}}
\date{Beijing Normal University}

\maketitle

\bigskip\bigskip

{\narrower{\narrower

\centerline{\bf Abstract}

\bigskip

 An
intrinsic multitype branching structure within the transient $(1,
R)$-RWRE is revealed.  The branching structure enables us to specify
the density of the absolutely continuous invariant measure for the
environments seen from the particle and reprove the LLN with an
drift {\it explicitly} in terms of the environment, comparing with
the results in Br\'emont (2002).

\bigskip


\bigskip


\par}\par}

\bigskip\bigskip


\section{Introduction and main results\label{s1}}

A random walk in random environment (RWRE, for short) $\{X_n, n\in
\mathbb{N}\}$ with bounded jumps on the line, written as
$(L,R)$-RWRE,  means that for each step the possible jump range to
the left is bounded by $L$ and to the right by $R$,  where $L$ and
$R$ are positive integers. The aim of this paper is  to reveal
 the branching structure within the $(1,R)$-RWRE,
 which is a  story different with but  essentially complement for that of
$(L, 1)$-RWRE (Hong and Wang, \cite{HW09}, 2009) when it  transient
to the right.

It is well-known that when the walk is transient, i.e., $X_n\to
\infty$, the intrinsic branching structure within the $(1,1)$-RWRE
is a Galton-Watson branching process with geometric offspring
distribution, which plays an important role in the proofs of the
limiting stable law (Kesten et al \cite{KKS75}, 1975), the renewal
theorem (Kesten \cite{K77}, 1977), and the law of large numbers
(Alili \cite{A99}, 1999) for this nearest random walk in random
environment, see also Gantert  and Shi (\cite{GS02}, 2002).

However, when $L$ or $R>1$, and the random walk is transient (we
consider $X_n\to \infty$),  It seems no result about the intrinsic
branching structure within the $(L, R)$-RWRE up to our knowledge,
even for random walks in non-random environment. Recently, partially
progress has been made by Hong and Wang (\cite{HW09}, 2009), where
 a multitype branching structure within the $(L,
1)$-RWRE have been revealed when $X_n\to \infty$ and the Kesten's
type stable law have been proved. In the present paper, we focus on
the opposite direction which is more complicated: we will figure out
the intrinsic branching structure within the $(1, R)$-RWRE when
$X_n\to \infty$(note that the situation for the $(L,1)$-RWRE
transient to the left is equivalent to that for the $(1,R)$-RWRE
transient to the right).
  We will discuss
$R=2$ in detail and extend to general $R$ at the end.

RWRE has been studied extensively  in recent years, especially for
the nearest $(1,1)$-RWRE on which many results have been obtained,
see for example \cite{S75}, \cite{KKS75}, \cite{K77}, \cite{A88},
etc.,  we refer to Sznitman (\cite{S02}, 2002) and  Zeitouni
(\cite{Z04}, 2004) as a general review. For the $(L, R)$-RWRE, Key
(\cite{K84}, 1984) discussed recurrence/transience criterion in
terms of the sign of two intermediate Liapounov exponents of a
random matrix. Br\'emont (\cite{B09}, 2009) formulate a criterion
for the existence of the absolutely continuous invariant measure for
the environments seen from the particle and deduce a
characterization of the non-zero-speed regime of the model. It
should note that the $(L, R)$-RWRE can be treated as a special case
of the random walks in random environments on a strip, see
Bolthausen  and  Goldsheid (\cite{BG00}, 2000; \cite{BG08}, 2008 )
and Roitershtein (\cite{R08}, 2008).

Br\'emont (\cite{B02}, 2002) has proved a recurrence/transience
criterion for $(L, 1)$-RWRE involving the greatest Lyapunov exponent
with respect to a random matrix $M$, and the law of large numbers by
assuming the $(IM)$ condition related the existence of an invariant
measure for the environments seen from the particle. As an important
application of our branching structure within the $(1, R)$-RWRE, it
enables us to specify the density of the absolutely continuous
invariant measure {\it explicitly}  and reprove the LLN with an
drift {\it explicitly} in terms of the environment.

We now introduce the model of $(1,2)$-RWRE. Generally speaking,
random walks in random environment involves two kinds of randomness:
the transition probability (we call it the ``environment"), which is
chosen from a specified distribution; and the random walk driven by
the chosen ``environment". Specifically, let $\Lambda=\{-1,1,2\}$ be
the set of possible jump range of the random walks. Let $\mathcal
{M}(\Lambda)$ be the collection of all probability measures on
$\Lambda$. Then define an environment to be an element
$\omega=\{(q(\omega)_x,p_1(\omega)_x,p_2(\omega)_x):x\in
\mathbb{Z}\}\in \mathcal {M}(\Lambda)^{\mathbb{Z}}=:\Omega$. Let $P$
be a stationary and ergodic probability distribution on $(\Omega,
\mathcal {F})$ and $\theta$ be the spatial shift, i.e.,
$(\theta\omega)_n=\omega_{n+1}$. Assume that $\theta$ is an
invertible transformation on the probability space
$(\Omega,\mathcal{F},P)$, measurable as well as its inverse and
preserving $P$. Moreover, assume that the environment is elliptic:
\beqnn
\exists ~\varepsilon>0,~ \forall~ z\in\{1,2\},~
(p_z/q)\geq\varepsilon,~ P\mbox{\textrm{-a.s.}}.
\eeqnn
Given an environment $\omega\in\Omega$, one can define a random walk
$\{X_n\}$ in the environment $\omega$ to be a time-homogeneous
Markov chain on $\mathbb{Z}$ with $X_0=0$ and the transition
probabilities
\beqnn
& &P_\omega(X_{n+1}=x-1|X_n=x):=q(\theta^x\omega)_0=q(\omega)_x,\\
& &P_\omega(X_{n+1}=x+z|X_n=x):=p_z(\theta^x\omega)_0=p_z(\omega)_x,
\eeqnn
for all $x\in \mathbb{Z}$ and $z\in \{1,2\}$. For each
$\omega$, we use $P_{\omega}$ to denote the law induced on the space
of paths $(\mathbb{Z}^\mathbb{N}, \mathcal{G})$. Then, define a
probability measure $\mathbb{P}:=P\otimes P_{\omega}$ on
$(\Omega\times\mathbb{Z}^\mathbb{N}, \mathcal{F}\times\mathcal{G})$
by
\beqnn
\mathbb{P}(F\times G):=\int_FP_{\omega}(G)P(\mbox{d}\omega),\quad
F\in\mathcal{F},~G\in\mathcal{G}.
\eeqnn
Statements involving
$P_{\omega}$ and $\mathbb{P}$ are called {\it quenched} and {\it
annealed}, respectively. Generally, with a slight abuse of notation,
$\mathbb{P}$ can also be used to denote the marginal on
$\mathbb{Z}^\mathbb{N}$. Expectations under $P_{\omega}$ and
$\mathbb{P}$ will be denoted by $E_\omega$ and $\mathbb{E}$,
respectively.

We will need the following notations. $p_z(\theta^x\omega)_0$ will
be simply denoted by $p_z(x)$. Any expression of the form
$f(\theta^k\omega)$ will be simply denoted by $f(k)$. For the random
walk $\{X_n\}$ and $k\in\mathbb{Z} \backslash\{0\}$, we write
$P_\omega^k$ for the quenched probability starting at $k$, and
$E_\omega^k$ for the corresponding expectation.

\

Assume that $X_0=0$ and the $(1,2)$-random walk is transient to the
right, i.e., $X_n\to\infty$, $\mathbb{P}$-a.s.. Let $T_0=0$, and
\beqnn T_k=\inf\{n> T_{k-1}:X_n>X_{T_{k-1}}\}, \quad k\ge 1 \eeqnn
be the sequence of {\it ladder times} of the random walks. Note that
$T_k<\infty,~\mathbb{P}$-a.s.. To calculate $T_1$ accurately, we
will figure out a multitype branching processes by decomposing the
path of the walk. Intuitively,  if the walk from $i\le0$ take a step
to $i-1$, it must crossing back to $i$ or jumping over $i$ (to
$i+1$) because of $T_k<\infty,~\mathbb{P}$-a.s., in which there are
only three kind of backing ways: from $i-1$ to $i$, from $i-2$ to
$i$ and from $i-1$ to $i+1$. So we divide all the steps   from $i $
to $i-1$ into three kind of steps according the crossing back ways.
Let $A(i),B(i)$, and $C(i)$ are the numbers of steps from $i$ to
$i-1$ before time $T_1$ with crossing-back from $i-1$ to $i$, $i-2$
to $i$ and $i-1$ to $i+1$, respectively. And for the last step of
$T_1$, we can consider it as a immigration for the multitype
branching processes.

Set for $i\le0$,
$$
U(i)=[A(i),B(i),C(i)],
$$
 Then we have the branching structure
within the $(1,2)$-RWRE as follows.

\btheorem\label{t1.1}
Assume $X_n\to \infty,~\mathbb{P}$-a.s.. Then for $P$-a.s. $\omega$,
$\Big(U(i)=[A(i),B(i),C(i)]\Big)_{i\le0}$ is an inhomogeneous
multitype branching process with immigration
\beqnn
U(1)=[1,0,0],~~\mbox{with probability}~~\frac{p_1(0)}{1-\alpha(0)-\beta(0)},\\
U(1)=[0,1,0],~~\mbox{with probability}~~\frac{\gamma(0)}{1-\alpha(0)-\beta(0)},\\
U(1)=[0,0,1],~~\mbox{with probability}~~\frac{p_2(0)}{1-\alpha(0)-\beta(0)}.
\eeqnn
The offspring distribution is given by
\beqlb
& &P_\omega\Big(U(i)=[a,b,0]~\Big|~U(i+1)=[1,0,0]\Big)=[1-\alpha(i)-\beta(i)]C_{a+b}^a\alpha(i)^a\beta(i)^b,\label{1.1}\\
& &P_\omega\Big(U(i)=[a,b,1]~\Big|~U(i+1)=[0,1,0]\Big)=[1-\alpha(i)-\beta(i)]C_{a+b}^a\alpha(i)^a\beta(i)^b,\label{1.2}\\
& &P_\omega\Big(U(i)=[a,b,0]~\Big|~U(i+1)=[0,0,1]\Big)=[1-\alpha(i)-\beta(i)]C_{a+b}^a\alpha(i)^a\beta(i)^b,\label{1.3}
\eeqlb
where
\beqnn
& &\gamma(i)=q(i)\cdot P_\omega^{i-1}[(-\infty,i-1),i+1],\\
& &\alpha(i)=q(i)\cdot P_\omega^{i-1}[(-\infty,i-1),i]\cdot\frac{p_1(i-1)}{p_1(i-1)+\gamma(i-1)},\\
& &\beta(i)=q(i)\cdot
P_\omega^{i-1}[(-\infty,i-1),i]\cdot\frac{\gamma(i-1)}{p_1(i-1)+\gamma(i-1)},
\eeqnn and $P_\omega^{i-1}[(-\infty,i-1),i-1+j]:= P_\omega^{i-1}
\{\mbox{reach $[i,+\infty)$ for the first time at the point
$i-1+j$}\}$ for $j=1,2$, the exit probabilities which can be
expressed in terms of the environment in Lemma \ref{l2.1}. \qed
\etheorem

As  an immediate consequence of Theorem \ref{t1.1}, We can get the
offspring quenched mean matrix of the  multitype branching process
and the quenched mean of the $T_1$ .

\bcorollary\label{c1.1}
Assume $X_n\to \infty,~\mathbb{P}$-a.s.. Then for $P$-a.s. $\omega$, and $i\le0$,
the offspring mean matrix of the $(-i+1)$-th generation
of the multitype branching process is
\beqnn
N(i)= \left(
\begin{array}{ccc}
\frac{\alpha(i)}{1-\alpha(i)-\beta(i)} & \frac{\beta(i)}{1-\alpha(i)-\beta(i)} &0\\
\frac{\alpha(i)}{1-\alpha(i)-\beta(i)} & \frac{\beta(i)}{1-\alpha(i)-\beta(i)} &1\\
\frac{\alpha(i)}{1-\alpha(i)-\beta(i)} & \frac{\beta(i)}{1-\alpha(i)-\beta(i)} &0
\end{array} \right).
\eeqnn
Moreover,
$$
E_\omega(T_1)=1+\Big\langle(2,2,1),~\frac{1}{1-\alpha(0)-\beta(0)}\Big(p_1(0),~\gamma(0),~p_2(0)\Big)
\cdot \sum_{i\le0}N(0)\cdots N(i)\Big\rangle.
$$
\qed
\ecorollary

As an application of the branching structure to the random walks in
non-random environment, let
$$
X_0=0,\quad X_n=\xi_1+\cdots+\xi_n,
$$
where $\xi_1,\ \xi_2,\cdots$ is a series of
i.i.d random variables with
\beqnn
P(\xi_1=-1)=q,\quad P(\xi_1=1)=p_1,\quad P(\xi_1=2)=p_2.
\eeqnn
The computability of $E(T_1)$ by the branching structure as in
Corollary \ref{c1.1} enable us to validate the Wald's
equality.

\bproposition\label{p1.1}
Assume that $E(\xi_1)=p_1+2p_2-q>0$. Then the Wald's equality holds:
$$E(X_{T_1})=E(T_1)\cdot E(X_1).$$
\qed
\eproposition

From the point of view ``environment viewed from the particles",
define $\overline{\omega}(n)=\theta^{X_n}\omega $. Then the process
$\{\overline{\omega}(n)\}$ is a Markov process under either
$P_\omega$ or $\mathbb{P}$, with state space $\Omega$ and transition
kernel
\beqnn
K(\omega, \mbox{d}\omega')=q{\bf1}_{\omega'=\theta^{-1}\omega}+p_1{\bf1}_{\omega'=\theta\omega}
+p_2{\bf1}_{\omega'=\theta^2\omega}.
\eeqnn
Set $\varphi^1_{\theta^k\omega}=P_{\theta^k\omega}(X_{T_1}=1)$, and
$\varphi^2_{\theta^k\omega}=P_{\theta^k\omega}(X_{T_1}=2)=1-\varphi_{\theta^k\omega}^1$.
Whenever $\mathbb{E}(T_1)<\infty$, define the measure
\beqnn
Q(B)=\mathbb{E}\bigg(\frac{{\bf1}_{X_{T_1}=1}}{\varphi_\omega^1}
\sum_{i=0}^{T_1-1}{\bf1}_{\{\overline{\omega}(i)\in B\}}+
\frac{{\bf1}_{X_{T_1}=2}}{\varphi_\omega^2}\sum_{i=0}^{T_1-1}{\bf1}_{\{\overline{\omega}(i)\in
B\}}\bigg),\quad \overline{Q}(B)= \frac{Q(B)}{Q(\Omega)}.
\eeqnn
The significance of the branching structure to the $(1, 2)$-RWRE is
that we can express the density of $Q$ with respect to $P$ {\it explicitly}.

\btheorem\label{t1.2}
Assume $\mathbb{E}(T_1)<\infty$. Then $Q(\cdot)$ is invariant under the Markov kernel $K$,
that is,
\beqnn
Q(B)=\iint{\bf1}_{\omega'\in B}K(\omega,\mbox{d}\omega')Q(\mbox{d}\omega).
\eeqnn
Furthermore,
\beqnn
\frac{\mbox{d}Q}{\mbox{d}P}=\Pi(\omega),
\eeqnn
where
\beqnn
\Pi(\omega)=\frac{2+\bigg\langle(1,1,1),
\sum_{i\ge1}\Big(\frac{p_1(i)}{p_1(i)+\gamma(i)},\frac{\gamma(i)}{p_1(i)+\gamma(i)},1\Big)\cdot
N(i)\cdots N(1)\bigg\rangle}{1-\alpha(0)-\beta(0)}.
\eeqnn
\qed
\etheorem

Therefore, the LLN for the $(1, 2)$-RWRE can be reproved with an {\it explicit} drift.

\btheorem\label{t1.3}
Assume $\mathbb{E}(T_1)<\infty$. Then
$$
\lim_{n\to\infty}\frac{X_n}{n}= v_P, \ \ ~\mathbb{P}\mbox{-a.s.},
$$
where
\beqnn
v_P=\frac{E_P\bigg[\frac{p_1(0)+2p_2(0)-p_{-1}(0)}{1-\alpha(0)-\beta(0)}
\Big(2+\Big\langle(1,1,1),\sum_{i\ge1}\Big(\frac{p_1(i)}{p_1(i)+\gamma(i)},
\frac{\gamma(i)}{p_1(i)+\gamma(i)},1\Big)\cdot N(i)\cdots N(1)\Big\rangle\Big)\bigg]}
{E_P\bigg[2+\Big\langle (2,2,1), \sum_{i\ge0}
\Big(\frac{p_1(i)}{p_1(i)+\gamma(i)},\frac{\gamma(i)}{p_1(i)+\gamma(i)},1\Big)\cdot
N(i)\cdots N(0)\Big\rangle \bigg]}.
\eeqnn \qed
\etheorem

\

We arrange the remainder of this paper as follows. In section
\ref{s2}, we formulate the intrinsic multitype branching structure
within the $(1, 2)$-RWRE under the assumption $X_n\to\infty$, i.e.,
Theorem \ref{t1.1}, Corollary \ref{c1.1} and Proposition \ref{p1.1}
are proved. In section \ref{s3}, we specify the density of the
absolutely continuous invariant measure for the environments seen
from the particle {\it explicitly} based on the branching structure,
i.e., Theorem \ref{t1.2} is proved. In section \ref{s4}, as an
application of the branching structure, the LLN will be reproved
with an {\it explicit} drift, that is, Theorem \ref{t1.3} is proved.
Finally in section \ref{s5}, for general $R>1$, we give the
intrinsic multitype branching structure within the $(1, R)$-RWRE
under the assumption $X_n\to\infty$.


\bigskip

\section{Intrinsic multitype branching structure within
the transient $(1,2)$-RWRE\label{s2}}

\setcounter{equation}{0}

We introduce the following exit probabilities of leaving a given
interval from the right side. Consider integers $a,b,k$ with $a\le
k\le b$, define
\beqnn
P_\omega^k[(a,b),b+1]= P_\omega^k \{\mbox{reach $(b,+\infty)$ before $(-\infty, a)$ and at the point $b+1$}\},\\
P_\omega^k[(a,b),b+2]= P_\omega^k \{\mbox{reach $(b,+\infty)$ before $(-\infty, a)$ and at the point $b+2$}\}.
\eeqnn
These exit probabilities play an important role in the offspring distribution
of the branching structure, which can be expressed in terms of the
environment (Br\'emont \cite{B02}, page.1273-4, Lemma 2.1 and
Proposition 2.2). Only a slight modification should be made: exiting
from the left in \cite{B02} corresponds exiting from the right here.
We still give the details of the proof for convenience.

\blemma\label{l2.1} For $n\ge 2$, we have \beqnn
P_\omega^{i}[(-n,i),i+1]= \frac{\langle e_1,[M(i)+\cdots+M(-n)\cdots
M(i)]v\rangle}
{1+\langle e_1, [M(i)+\cdots+M(-n)\cdots M(i)]e_1\rangle},\\
P_\omega^{i}[(-n,i),i+2]= \frac{\langle e_1,[M(i)+\cdots+M(-n)\cdots
M(i)]e_2\rangle} {1+\langle e_1,[M(i)+\cdots+M(-n)\cdots
M(i)]e_1\rangle}, \eeqnn where $e_1=(1,0)'$, $e_2=(0,1)'$,
$v=e_1-e_2$ and \beqnn M(i):= \left(
\begin{array}{cc}
\frac{p_1(i)+p_2(i)}{q(i)} & \frac{p_2(i)}{q(i)}\\
1 & 0
\end{array} \right).
\eeqnn Furthermore, if $X_n\to+\infty,~\mathbb{P}$-a.s., then
\beqlb\label{2.1}
P_\omega^{i}[(-\infty,i),i+1]+P_\omega^{i}[(-\infty,i),i+2]=1.
\eeqlb \elemma

{\it Proof}. Set $f(k)=P_\omega^k[(-n,i),i+1]$. For $k$ such that
$-n\le k\le i$, using Markov property, we get a harmonic-type
recurrence equation:
$$f(k)=p_1(k)f(k+1)+p_2(k)f(k+2)+q(k)f(k-1),$$
with $f(i+1)=1,~f(k)=0$ for $k\ge i+2$ or $k\le-n-1$. Setting
$g(k)=f(k)-f(k-1)$, we obtain \beqnn
g(k)&=&\frac{p_1(k)}{q(k)}g(k+1)+\frac{p_2(k)}{q(k)}\Big(g(k+2)+g(k+1)\Big)\\
&=&\frac{p_1(k)+p_2(k)}{q(k)}g(k+1)+\frac{p_2(k)}{q(k)}g(k+2).
\eeqnn
Rewrite the equation as
\beqnn
\left(\begin{array}{c}
g(k)\\
g(k+1)
\end{array}\right)
=\left( \begin{array}{cc}
\frac{p_1(k)+p_2(k)}{q(k)}&\frac{p_2(k)}{q(k)}\\
1 & 0
\end{array} \right)
\left(\begin{array}{c}
g(k+1) \\
g(k+2)
\end{array}\right).
\eeqnn
Set
\beqnn
U(k)= \left( \begin{array}{c}
g(k)\\
g(k+1)
\end{array}\right),
\quad \mbox{and}\quad
M(k) = \left( \begin{array}{cc}
a_1(k) & a_2(k) \\
1  &0
\end{array} \right),
\eeqnn
where $a_1(k)=\frac{p_1(k)+p_2(k)}{q(k)}$,
$a_2(k)=\frac{p_2(k)}{q(k)}$. Thus, we obtain the relations
\beqnn
U(k)&=&M(k)U(k+1)\\
&=&M(k)\cdots M(i)U(i+1),
\eeqnn
with $U(i+1)=([1-f(i)],-1)'$ and $U(-n)=(f(-n),[f(-n+1)-f(-n)])'$.
Summing from $-n$ to $i$, we deduce that
\beqnn
f(i)&=&\Big\langle e_1, M(i)U(i+1)\Big\rangle+\cdots
+\Big\langle e_1, M(-n)M(-n+1)\cdots M(i)U(i+1)\Big\rangle\\
&=&\bigg\langle e_1, \Big[M(i)+\cdots +M(-n)\cdots M(i)\Big]\cdot
\left(
\begin{array}{c}
1-f(i)\\
-1
\end{array}\right)\bigg\rangle.
\eeqnn
The first formula then follows. By similar reasoning, one can get
the second formula. The second part of the lemma is just the conclusion
of Proposition 2.2 in \cite{B02}. \qed

\bigskip

Now we introduce the branching structure. Recall $X_0=0$ and the
sequence of {\it ladder times} of the random walks:
$T_0=0$ and
\beqnn
T_k=\inf\{n> T_{k-1}:X_n>X_{T_{k-1}}\},\quad k\ge1.
\eeqnn
Note that $T_k<\infty,~\mathbb{P}$-a.s. if $X_n\to+\infty,~\mathbb{P}$-a.s..
Define, for $ i\le0$,
\beqnn
& &\eta_{i,0}=\min\{k\le T_1:~X_k=i\},\\
& &\theta_{i,0}=\min\{\eta_{i,0}<k\le T_1:~X_{k-1}=i,X_k=i-1\},
\eeqnn
and for $j\ge1$,
\beqnn
& &\alpha_{i,j}=\min\{\theta_{i,j-1}<k\le T_1:~X_{k-1}=i-1, X_k=i\},\\
& &\beta_{i,j}=\min\{\theta_{i,j-1}<k\le T_1:~X_{k-1}=i-2,X_k=i\},\\
& &\gamma_{i,j}=\min\{\theta_{i,j-1}<k\le T_1:~X_{k-1}=i-1,X_k=i+1\},\\
& &\eta_{i,j}=\min\{\alpha_{i,j},\beta_{i,j},\gamma_{i,j}\},\\
& &\theta_{i,j}=\min\{\eta_{i,j}<k\le T_1:~X_{k-1}=i,X_k=i-1\},
\eeqnn
with the usual convention that the minimum over an empty set
is $+\infty$. We refer to the time interval
$[\theta_{i,j-1},\eta_{i,j}]$ as the $j$-th excursion from $i-1$ to
$\{i,i+1\}$. For any $j\ge0$, any $i\le0$, define
\beqnn
& &A_{i,j}=\#\{k\ge0: \theta_{i+1,j}<\theta_{i,k}<\eta_{i+1,j+1}, ~\mbox{and}~\eta_{i,k+1}=\alpha_{i,k+1}\},\\
& &B_{i,j}=\#\{l\ge0: \theta_{i+1,j}<\theta_{i,l}<\eta_{i+1,j+1}, ~\mbox{and}~\eta_{i,l+1}=\beta_{i,l+1}\},\\
& &C_{i,j}=\#\{m\ge0: \theta_{i+1,j}<\theta_{i,m}<\eta_{i+1,j+1},
~\mbox{and}~\eta_{i,m+1}=\gamma_{i,m+1}\}.
\eeqnn
Note that $A_{i,j}, B_{i,j}$, and $C_{i,j}$ are the numbers of steps from $i$
to $i-1$ during the $(j+1)$-th excursion from $i$ to $\{i+1,i+2\}$
with crossing-back from $i-1$ to $i$, $i-2$ to $i$ and $i-1$ to
$i+1$, respectively. Special attentions should be paid to the
different ending ways of each excursion and their consequences, take
the $(j+1)$-th excursion from $i$ to $\{i+1,i+2\}$ for example:
\begin{itemize}
   \item If $\eta_{i+1,j+1}=\alpha_{i+1,j+1}$, then the excursion
   $[\theta_{i+1,j+1},\alpha_{i+1,j+1}]$ ends by jumping from $i$ to $i+1$.
   \item If $\eta_{i+1,j+1}=\gamma_{i+1,j+1}$, then the excursion
   $[\theta_{i+1,j+1},\gamma_{i+1,j+1}]$ ends by jumping from $i$ to $i+2$.
\end{itemize}
Furthermore, in the above two cases, excursions from $i-1$ to
$\{i,i+1\}$ included in the $(j+1)$-th excursion from $i$ to
$\{i+1,i+2\}$ must end at $i$ and consequently, $C_{i,j}=0$.
\begin{itemize}
   \item If $\eta_{i+1,j+1}=\beta_{i+1,j+1}$, then the excursion
   $[\theta_{i+1,j+1},\beta_{i+1,j+1}]$ ends by jumping from $i-1$ to $i+1$.
\end{itemize}
Moreover, in this case, the last excursion from $i-1$ to $\{i,i+1\}$
included in the $j+1$-th excursion from $i$ to $\{i+1,i+2\}$ also
ends by jumping from $i-1$ to $i+1$, that is, if $\theta_{i,m_0}$ is
the beginning of it, then
$\eta_{i,m_0+1}=\gamma_{i,m_0+1}=\beta_{i+1,j+1}=\eta_{i+1,j+1}$.
Therefore, $C_{i,j}=1$.

\

Define for $i\le0$,
$$
A(i)=\sum_{j\ge0}A_{i,j},~B(i)=\sum_{j\ge0}B_{i,j},~C(i)=\sum_{j\ge0}C_{i,j}
$$
to be the numbers of steps from $i$ to $i-1$ before time $T_1$ with
crossing-back from $i-1$ to $i$, $i-2$ to $i$ and $i-1$ to $i+1$,
respectively. Define for $i\le0$,
$$
U(i)=[A(i),B(i),C(i)].
$$
Then $A(i)+B(i)+C(i)$ is the total
number of steps the walk jumping from $i$ to $i-1$ before $T_1$. The
number of steps crossing back from $i-1$ to $i$ is $A(i)+B(i)$, because
$X_n\to \infty,~\mathbb{P}$-a.s. and $C(i)=B(i+1)$ is counted as steps
crossing back from $i$ to $i+1$. Thus we have
\beqlb\label{2.2}
T_1&=&1+\sum_{i\le0}\Big(2A(i)+2B(i)+C(i)\Big)\nonumber \\
&=&1+\Big\langle(2,2,1), \sum_{i\le0}U(i)\Big\rangle.
\eeqlb
In order to study $T_1$, we consider $\{U(i)\}_{i\le0}$ instead. We
need the following probabilities. For $i\le0$,
\beqlb
\alpha(i)&=&P^i_\omega[\mbox{ the walk jumping from $i$ to $i-1$, and finally jumping back }\label{2.3}\\
& &\mbox{from $i-1$ to $i$ before $T_1$}],\nonumber\\
\beta(i)&=&P^i_\omega[\mbox{ the walk jumping from $i$ to $i-1$, and finally jumping back }\label{2.4}\\
& &\mbox{from $i-2$ to $i$ before $T_1$}],\nonumber\\
\gamma(i)&=&P^i_\omega[\mbox{ the walk jumping from $i$ to $i-1$, and finally jumping back }\label{2.5}\\
& &\mbox{from $i-1$ to $i+1$ before $T_1$}].\nonumber
\eeqlb
Expressions of these probabilities will be given later in terms of the exit
probabilities. By Markov property, and the different ending cases
explained above, we obtain, for $i\le0$,
\beqlb
& &P_\omega\Big(U(i)=[a,b,0],U(i+1)=[1,0,0]~\Big|~|U(i+1)|=1\Big)=C_{a+b}^a\alpha(i)^a\beta(i)^bp_1(i),\label{2.6}\\
& &P_\omega\Big(U(i)=[a,b,1],U(i+1)=[0,1,0]~\Big|~|U(i+1)|=1\Big)=C_{a+b}^a\alpha(i)^a\beta(i)^b\gamma(i),\label{2.7}\\
& &P_\omega\Big(U(i)=[a,b,0],U(i+1)=[0,0,1]~\Big|~|U(i+1)|=1\Big)=C_{a+b}^a\alpha(i)^a\beta(i)^bp_2(i).\label{2.8}
\eeqlb
Additionally, for $i=0$, we define $U(1)=[A(1),B(1),C(1)]$
for consistency. Noting that the walk starts at $0$, and the
excursion from $0$ to $\{1,2\}$ ends at time $T_1$, there is only
one jump crossing up from $0$ to $1$ in the time interval $[0,T_1]$,
which can be regarded as $|U(1)|=1$ if $T_1<\infty,~\mathbb{P}$-a.s..
If the excursion from $0$ to $\{1,2\}$ ends by jumping from $0$ to
$1$, then set $A(1)=1$, in other words, $U(1)=[1,0,0]$. By similar
reasoning, $U(1)=[0,1,0]$ if the ending jump is from $-1$ to $1$ and
$U(1)=[0,0,1]$ if the ending jump is from $0$ to $2$. Summing over
$a,b\ge0$ in (\ref{2.6}), we have
\beqlb\label{2.9}
\nonumber & &  P_\omega\Big(U(i+1)=[1,0,0]~\Big|~|U(i+1)|=1\Big)\\
\nonumber &=&\sum_{a,b\ge0} P_\omega\Big(U(i)=[a,b,0], U(i+1)=[1,0,0]~\Big|~|U(i+1)|=1\Big)\\
&=&\sum_{a,b\ge0}C_{a+b}^a\alpha(i)^a\beta(i)^bp_1(i)=\frac{p_1(i)}{1-\alpha(i)-\beta(i)}.
\eeqlb
Similarly,
\beqlb
& &P_\omega\Big(U(i+1)=[0,1,0]~\Big|~|U(i+1)|=1\Big)=\frac{\gamma(i)}{1-\alpha(i)-\beta(i)},\label{2.10}\\
& &P_\omega\Big(U(i+1)=[0,0,1]~\Big|~|U(i+1)|=1\Big)=\frac{p_2(i)}{1-\alpha(i)-\beta(i)}.\label{2.11}
\eeqlb

Now we are ready to calculate $\alpha(i),~\beta(i)$, and
$\gamma(i)$. Firstly, by the definitions (\ref{2.3})--(\ref{2.5}),
we can see
\beqnn
& &\alpha(i)=P_\omega\Big(U(i)=[1,0,0]~\Big|~|U(i)|=1\Big),\\
& &\beta(i)=P_\omega\Big(U(i)=[0,1,0]~\Big|~|U(i)|=1\Big),\\
& &\gamma(i)=P_\omega\Big(U(i)=[0,0,1]~\Big|~|U(i)|=1\Big).
\eeqnn
Thus from (\ref{2.9}) and (\ref{2.10}), we have the ratio
\beqlb\label{2.12}
\alpha(i):\beta(i)=p_1(i-1):\gamma(i-1).
\eeqlb
On the other hand, by the definitions (\ref{2.3})--(\ref{2.5})
again, we know
\beqnn
\alpha(i)+\beta(i)=P^i_\omega[\mbox{ the walk jumping from $i$ to $i-1$,}~~~~~~~~~~~~~~~~~~~~\\
\mbox{ and finally jumping back to $i$ before $T_1$}],\\
\gamma(i)=P^i_\omega[\mbox{ the walk jumping from $i$ to $i-1$,}~~~~~~~~~~~~~~~~~~~~\\
\mbox{ and finally jumping back to $i+1$ before $T_1$}].
\eeqnn
Recall the exit probabilities, we obtain
\beqlb
\alpha(i)+\beta(i)&=&q(i)\cdot P_\omega^{i-1}[(-\infty,i-1),i],\label{2.13}\\
\gamma(i)&=&q(i)\cdot P_\omega^{i-1}[(-\infty,i-1),i+1].\label{2.14}
\eeqlb

Finally, combining (\ref{2.12}) and (\ref{2.13}) we get
\beqlb
& &\alpha(i)=q(i)\cdot P_\omega^{i-1}[(-\infty,i-1),i]\cdot\frac{p_1(i-1)}{p_1(i-1)+\gamma(i-1)},\label{2.15}\\
& &\beta(i)=q(i)\cdot P_\omega^{i-1}[(-\infty,i-1),i]\cdot\frac{\gamma(i-1)}{p_1(i-1)+\gamma(i-1)}.\label{2.16}
\eeqlb

\

\noindent{\it Proof of Theorem \ref{t1.1}.} By
(\ref{2.6})--(\ref{2.8}) and (\ref{2.9})--(\ref{2.11}), we obtain
\beqlb
& &P_\omega\Big(U(i)=[a,b,0]~\Big|~U(i+1)=[1,0,0]\Big)=[1-\alpha(i)-\beta(i)]C_{a+b}^a\alpha(i)^a\beta(i)^b,\label{2.17}\\
& &P_\omega\Big(U(i)=[a,b,1]~\Big|~U(i+1)=[0,1,0]\Big)=[1-\alpha(i)-\beta(i)]C_{a+b}^a\alpha(i)^a\beta(i)^b,\label{2.18}\\
& &P_\omega\Big(U(i)=[a,b,0]~\Big|~U(i+1)=[0,0,1]\Big)=[1-\alpha(i)-\beta(i)]C_{a+b}^a\alpha(i)^a\beta(i)^b.\label{2.19}
\eeqlb
(\ref{2.17})--(\ref{2.19}) give the offspring distribution for the
$(-i+1)$-th generation ($i\leq 0$). Indeed, it is enough to show
\beqnn
& &\sum_{a,b\ge0}\bigg[P_\omega\Big(U(i)=[a,b,0]~\Big|~U(i+1)=[1,0,0]\Big)\cdot P_\omega\Big(U(i+1)=[1,0,0]~\Big|~|U(i+1)|=1\Big)\\
& &\qquad +P_\omega\Big(U(i)=[a,b,1]~\Big|~U(i+1)=[0,1,0]\Big)\cdot P_\omega\Big(U(i+1)=[0,1,0]~\Big|~|U(i+1)|=1\Big)\\
& &\qquad +P_\omega\Big(U(i)=[a,b,0]~\Big|~U(i+1)=[0,0,1]\Big)\cdot
P_\omega\Big(U(i+1)=[0,0,1]~\Big|~|U(i+1)|=1\Big)\bigg]=1,
\eeqnn
which is equivalent to check
$$
\frac{p_1(i)+\gamma(i)+p_2(i)}{1-\alpha(i)-\beta(i)}=1.
$$
By (\ref{2.1}), (\ref{2.13}), and (\ref{2.14}), we know
$\alpha(i)+\beta(i)+\gamma(i)=q(i)$. The conclusion therefore
follows from $p_1(i)+p_2(i)+q(i)=1$.

Since $X_0=0$ and $X_n\to\infty,~\mathbb{P}$-a.s., we have
$T_1<+\infty,~\mathbb{P}$-a.s., so we imagine that there are different
types of particles immigrating into the system $U(1)=[A(1),B(1),C(1)]$
and hence $|U(1)|=1~\mathbb{P}$-a.s..
By (\ref{2.9}), we know
\beqnn
P_\omega\Big(U(1)=[1,0,0]\Big)=P_\omega\Big(U(1)=[1,0,0]~\Big|~|U(1)|=1\Big)=\frac{p_1(0)}{1-\alpha(0)-\beta(0)}.
\eeqnn
By similar argument,
\beqnn
& &P_\omega\Big(U(1)=[0,1,0]\Big)=\frac{\gamma(0)}{1-\alpha(0)-\beta(0)},\\
& &P_\omega\Big(U(1)=[0,0,1]\Big)=\frac{p_2(0)}{1-\alpha(0)-\beta(0)}.
\eeqnn
Thus the immigration of the multi-type branching process
follows.  This completes the proof.   \qed

\

\noindent{\it Proof of Corollary  \ref{c1.1}.} Summing over $b\ge0$
in (\ref{2.17}), we obtain
\beqnn
P_\omega\Big(A(i)=a~\Big|~U(i+1)=[1,0,0]\Big)
=\frac{1-\alpha(i)-\beta(i)}{1-\beta(i)}\Big(\frac{\alpha(i)}{1-\beta(i)}\Big)^a.
\eeqnn
So the expected number of type-A offspring of a single type-A
particle in one generation is
\beqnn
N_{11}(i)=\sum_{a\ge0}a\cdot P_\omega \Big(A(i)=a~\Big|~U(i+1)=[1,0,0]\Big)
=\frac{\alpha(i)}{1-\alpha(i)-\beta(i)}.
\eeqnn
By the same argument, we can get all the elements of the
offspring mean matrix $N(i)$. For the second result, noting that by
the multitype branching process, we have
$$
E_\omega [U(i)]=E_\omega [U(i+1)]\cdot N(i)=E_\omega[U(1)]\cdot N(0)\cdots N(i),
$$
and
\beqnn
& &E_\omega[U(1)]\\
& &=[1,0,0]\cdot P_\omega\Big(U(1)=[1,0,0]\Big)+[0,1,0]\cdot P_\omega\Big(U(1)=[0,1,0]\Big)\\
& &\qquad\qquad\qquad\qquad\qquad\qquad\qquad\qquad\qquad\qquad+[0,0,1]\cdot P_\omega\Big(U(1)=[0,0,1]\Big)\\
& &=\Big(\frac{p_1(0)}{1-\alpha(0)-\beta(0)},
\frac{\gamma(0)}{1-\alpha(0)-\beta(0)},\frac{p_2}{1-\alpha(0)-\beta(0)}\Big)\\
& &=\frac{1}{1-\alpha(0)-\beta(0)}\Big(p_1(0),~\gamma(0),~p_2(0)\Big).
\eeqnn
The desired conclusion follows by taking the quenched expectation in (\ref{2.2}). \qed

\

As an immediate application of the branching structure, consider
random walks in non-random environment.  Let $\xi_1,\xi_2,\cdots$ be
independent variables with common distribution
\beqnn
P(\xi_1=-1)=q,\quad P(\xi_1=1)=p_1,\quad P(\xi_1=2)=p_2.
\eeqnn
The induced random walk is the sequence of random variables
$$
X_0=0,\quad X_n=\xi_1+\cdots+\xi_n.
$$

Define $T_1=\inf\{n:X_n>0\}$ as before. The countability of $T_1$ in
Corollary  \ref{c1.1} enable us to check the Wald's equality
(\cite{F71}, pg. 397) directly.

\bigskip

\noindent{\it Proof of Proposition \ref{p1.1}.}~ First, we have
$$
E(X_{T_1})=1\cdot P(X_{T_1}=1)+2\cdot P(X_{T_1}=2).
$$
By Lemma 2.1, we obtain
\beqlb
& &P(X_{T_1}=1)=P^0[(-\infty,0),1]=\lim_{n\to\infty} \frac{\langle e_1,
[M+M^2+\cdots+M^n]v\rangle} {1+\langle e_1, [M+M^2+\cdots+M^n]e_1\rangle},\label{2.20}\\
& &P(X_{T_1}=2)=P^0[(-\infty,0),2]=\lim_{n\to\infty} \frac{\langle e_1,
[M+M^2+\cdots+M^n]e_2\rangle} {1+\langle e_1, [M+M^2+\cdots+M^n]e_1\rangle},\label{2.21}
\eeqlb
where
\beqnn
M:=
\left(
\begin{array}{cc}
\frac{p_1+p_2}{q}&\frac{p_2}{q}\\
1&0
\end{array} \right).
\eeqnn
Note that
$$
M^n=A\Lambda^n A^{-1}
=\left(
\begin{array}{cc}
\lambda_1 & \lambda_2 \\
1 & 1 \\
\end{array}
\right)
\cdot\left(
\begin{array}{cc}
\lambda_1^n & 0 \\
0 & \lambda_2^n \\
\end{array}\right)
\cdot \frac{1}{\lambda_1-\lambda_2}\left(
\begin{array}{cc}
1& -\lambda_2 \\
-1 & \lambda_1 \\
\end{array}\right),
$$
where $\lambda_{1,2}=\frac{p_1+p_2\pm\sqrt{(p_1+p_2)^2+4qp_2}}{2q}$
are the eigenvalues of $M$. Then
\beqnn
& &\langle e_1, M^nv\rangle=(1+\lambda_2)\lambda_1^{n+1}-(1+\lambda_1)\lambda_2^{n+1},\\
& &\langle e_1, M^ne_2\rangle=(-\lambda_2)\lambda_1^{n+1}+\lambda_1\lambda_2^{n+1},\\
& &\langle e_1, M^ne_1\rangle=\lambda_1^{n+1}-\lambda_2^{n+1}.
\eeqnn
Since $EX_1>0$, we have $\lambda_1>1$, and
$\lambda_2\in(-1,0)$. Hence, by (\ref{2.20}) and (\ref{2.21}), we obtain
\beqnn
& &P(X_{T_1}=1)=1+\lambda_2,\\
& &P(X_{T_1}=2)=-\lambda_2.
\eeqnn
Therefore, $E(X_{T_1})=1-\lambda_2$.

\medskip

The next step is to calculate $E(T_1)$, which is done by the
branching process (Corollary \ref{c1.1}). Since the environment is
not random, there is no site index or integration with respect to
$P$ any more. Applying Corollary \ref{c1.1} to this random walk, we
have
$$
E(T_1)=1+\Big\langle (2,2,1),
\frac{1}{1-\alpha-\beta}\Big(p_1,~\gamma,~p_2\Big)\cdot \sum_{n\ge1}N^n\Big\rangle.
$$
where
\beqnn
N:= \left( \begin{array}{ccc}
\frac{\alpha}{1-\alpha-\beta} & \frac{\beta}{1-\alpha-\beta} &0\\
\frac{\alpha}{1-\alpha-\beta} & \frac{\beta}{1-\alpha-\beta} &1\\
\frac{\alpha}{1-\alpha-\beta} &
\frac{\beta}{1-\alpha-\beta} &0
\end{array} \right),
\eeqnn
Using the eigenvalues and eigenvectors of $N$, and noting
that the norm of the greatest eigenvalue of $N$ is less then $1$, we
have
\beqlb\label{2.22}
\sum_{n\ge1}N^n=\frac{1}{1-2\alpha-3\beta}
\begin{pmatrix}
\alpha & \beta & \beta\\
2\alpha & 2\beta & 1-2\alpha-\beta\\
\alpha & \beta &\beta
\end{pmatrix}.
\eeqlb
Thus,
$$
E(T_1)=\frac{\gamma+1-\alpha-\beta}{(1-\alpha-\beta)(1-2\alpha-3\beta)}.
$$
Substituting $\gamma=-q\lambda_2$,
$\alpha=q(1+\lambda_2)\frac{p_1}{p_1-q\lambda_2}$, and
$\beta=q(1+\lambda_2)\frac{-q\lambda_2}{p_1-q\lambda_2}$, we have
\beqnn
\gamma+1-\alpha-\beta&=&\Delta,\\
1-\alpha-\beta&=&\frac{1}{2}(1-q+\Delta),\\
1-2\alpha-3\beta&=&1-\frac{(1+q-\Delta)(p_1-3p_2+3\Delta)}{2(p_1-p_2+\Delta)},
\eeqnn
where $\Delta=\sqrt{(p_1+p_2)^2+4p_2q}$.
Then
\beqnn
E(T_1)
&=&\frac{2\Delta}{(1-q+\Delta)-\frac{(1-q+\Delta)(1+q-\Delta)(p_1-3p_2+3\Delta)}{2(p_1-p_2+\Delta)}}\\
&=&\frac{2\Delta}{(1-q+\Delta)-\frac{2q(p_1-p_2+\Delta)(p_1-3p_2+3\Delta)}{2(p_1-p_2+\Delta)}}\\
&=&\frac{2\Delta}{1-q-qp_1+3qp_2+(1-3q)\Delta}=:\frac{a(p_1,p_2,q)}{b(p_1,p_2,q)},
\eeqnn
and
$$
\frac{E(X_{T_1})}{E(X_1)}
=\frac{1-\lambda_2}{p_1+2p_2-q}=\frac{3q-1+\Delta}{2q(p_1+2p_2-q)}=:\frac{c(p_1,p_2,q)}{d(p_1,p_2,q)}.
$$
In order to prove $E(X_{T_1})=E(T_1)\cdot E(X_1)$, it suffices to show
$$
a(p_1,p_2,q)d(p_1,p_2,q)=b(p_1,p_2,q)c(p_1,p_2,q).
$$
In fact,
\beqnn
a(p_1,p_2,q)d(p_1,p_2,q)&=&4q(p_1+2p_2-q)\Delta\\
&=&4q(1-p_2-q+2p_2-q)\Delta\\
&=&4q(1+p_2-2q)\Delta.
\eeqnn
And,
\beqnn
& &b(p_1,p_2,q)c(p_1,p_2,q)\\
&=&[(1-q-qp_1+3qp_2)+(1-3q)\Delta](3q-1+\Delta)\\
&=&(1-q-qp_1+3qp_2)(3q-1)+(1-3q)\Delta^2+[(1-q-qp_1+3qp_2)+(1-3q)(3q-1)]\Delta\\
&=&(3q-1)(1-q-qp_1+3qp_2-(1-q)^2-4qp_2)\\
& &\qquad\qquad\qquad\qquad\qquad\qquad\qquad\qquad+[1-q-q(1-q-p_2)+3qp_2-(1-6q+9q^2)]\Delta\\
&=&(3q-1)(q(1-p_1-p_2)-q^2)+[4q+4qp_2-8q^2]\Delta\\
&=&4q(1+p_2-2q)\Delta.
\eeqnn
Thus, the desired conclusion follows. \qed

\bigskip

\noindent{\bf Remark}  Proposition \ref{p1.1} is a strong evidence
to validate the branching structure.


\bigskip

\section{Density of the absolutely continuous
invariant measure\label{s3}}

\setcounter{equation}{0}

Now we introduce the machinery of the ``environment viewed from
the particle". The first step consists of introducing an auxiliary
Markov chain. Starting from the RWRE $\{X_n\}$, define
$\overline{\omega}(n)=\theta^{X_n}\omega$. The sequence
$\{\overline{\omega}(n)\}$ is a process with paths in
$\Omega^\mathbb{N}$. This process is in fact a Markov process. The
proof is the same as Zeitouni (\cite{Z04}, page 204, Lemma 2.1.18),
we omit the details.

\blemma\label{l3.1}
The process $\{\overline{\omega}(n)\}$ is a
Markov process under either $P_\omega$ or $\mathbb{P}$, with
state space $\Omega$ and transition kernel
\beqnn
K(\omega, \mbox{d}\omega')=q{\bf1}_{\omega'=\theta^{-1}\omega}+p_1{\bf1}_{\omega'=\theta\omega}
+p_2{\bf1}_{\omega'=\theta^2\omega}.
\eeqnn
\qed
\elemma

The next step is to construct an invariant measure for the
transition kernel $K$. Assume that
$X_n\to\infty,~\mathbb{P}$-a.s., implying
$T_n<\infty,~\mathbb{P}$-a.s.. Set
$\varphi^1_{\theta^k\omega}=P_{\theta^k\omega}(X_{T_1}=1)$, and
$\varphi^2_{\theta^k\omega}=P_{\theta^k\omega}(X_{T_1}=2)=1-\varphi_{\theta^k\omega}^1$.
Whenever $\mathbb{E}(T_1)<\infty$, define the measure
\beqnn
Q(B)=\mathbb{E}\bigg(\frac{{\bf1}_{X_{T_1}=1}}{\varphi_\omega^1}\sum_{i=0}^{T_1-1}{\bf1}_{\{\overline{\omega}(i)\in
B\}}+\frac{{\bf1}_{X_{T_1}=2}}{\varphi_\omega^2}\sum_{i=0}^{T_1-1}{\bf1}_{\{\overline{\omega}(i)\in
B\}}\bigg),\quad \overline{Q}(B)= \frac{Q(B)}{Q(\Omega)}.
\eeqnn
Note that $\overline{Q}$ is a probability measure.

\

\noindent{\it Proof of Theorem \ref{t1.2}---invariant measure}.
We will show
\beqnn
Q(B)=\iint{\bf1}_{\omega'\in B}K(\omega,\mbox{d}\omega')Q(\mbox{d}\omega).
\eeqnn
On one hand,
\beqnn
Q(B)=\sum_{i=0}^\infty\mathbb{E}\bigg(\frac{{\bf1}_{X_{T_1}=1}}{\varphi_\omega^1};
T_1>i; \overline\omega(i)\in B\bigg)
    +\sum_{i=0}^\infty\mathbb{E}\bigg(\frac{{\bf1}_{X_{T_1}=2}}{\varphi_\omega^2};
T_1>i; \overline\omega(i)\in B\bigg).
\eeqnn
On the other hand,
\beqnn
& &\iint {\bf1}_{\omega'\in B}K(\omega,\mbox{d}\omega')Q(\mbox{d}\omega)\\
&=&\sum_{i=0}^\infty\mathbb{E}\bigg(\frac{{\bf1}_{X_{T_1}=1}}{\varphi_\omega^1};~ T_1>i;~
\overline\omega(i+1)\in B \bigg)
  +\sum_{i=0}^\infty\mathbb{E}\bigg(\frac{{\bf1}_{X_{T_1}=2}}{\varphi_\omega^2};~ T_1>i;~
\overline\omega(i+1)\in B \bigg)\\
&=&\sum_{j=1}^\infty\mathbb{E}\bigg(\frac{{\bf1}_{X_{T_1}=1}}{\varphi_\omega^1};~ T_1>j;~
\overline\omega(j)\in B \bigg)
  +\sum_{j=1}^\infty\mathbb{E}\bigg(\frac{{\bf1}_{X_{T_1}=2}}{\varphi_\omega^2};~ T_1>j;~
\overline\omega(j)\in B \bigg)\\
& &\qquad+\mathbb{E}\bigg(\frac{{\bf1}_{X_{T_1}=1}}{\varphi_\omega^1};~T_1<\infty;~
\overline\omega(T_1)\in B \bigg)
         +\mathbb{E}\bigg(\frac{{\bf1}_{X_{T_1}=2}}{\varphi_\omega^2};~T_1<\infty;~
\overline\omega(T_1)\in B \bigg).
\eeqnn
It only needs to show
\beqnn
 \mathbb{E}\bigg(\frac{{\bf1}_{X_{T_1}=i}}{\varphi_\omega^i};~ T_1>0;~ \overline\omega(0)\in B\bigg)
=\mathbb{E}\bigg(\frac{{\bf1}_{X_{T_1}=i}}{\varphi_\omega^i};~ T_1<\infty;~ \overline\omega(T_1)\in B\bigg),~i=1,2.
\eeqnn
Indeed,
\beqnn
& &\mbox{R.S.}=E_P\bigg[\frac{1}{\varphi_\omega^i}E_\omega\Big({\bf1}_{X_{T_1}=i},~ \overline\omega(T_1)\in B, T_1<\infty\Big)\bigg]\\
& &\qquad=E_P\bigg[P_\omega\Big(\overline\omega(T_1)\in B,~T_1<\infty\Big|X_{T_1}=i\Big)\bigg]\\
& &\qquad=~P~\bigg[\theta^i\omega\in B\cdot P_\omega\Big(T_1<\infty\Big|X_{T_1}=i\Big)\bigg]\\
& &\qquad=~P~\Big( \theta^i\omega\in B\Big)\mbox{ (since }P_\omega(T_1<\infty)=1\mbox{)}.
\eeqnn
Similarly, L.S.$=P(\omega\in B)$. Then, by the invariance of the
environment, we get L.S.=R.S.. This completes the proof.\qed

\

\noindent{\it Proof of Theorem \ref{t1.2}---density}. Let
$f:\Omega\to\mathbb{R}$ be measurable. Then,
\beqnn
\int f \mbox{d}Q&=&\mathbb{E}\bigg(\sum_{i=0}^{T_1-1}f(\overline\omega(i))
\frac{{\bf1}_{X_{T_1}=1}}{\varphi_\omega^1}+\sum_{i=0}^{T_1-1}f(\overline\omega(i))
\frac{{\bf1}_{X_{T_1}=2}}{\varphi_\omega^2}\bigg)\\
&=&\mathbb{E}\bigg(\sum_{i\le0}f(\theta^i\omega)V_i\frac{{\bf1}_{X_{T_1}=1}}{\varphi_\omega^1}
+\sum_{i\le0}f(\theta^i\omega)V_i\frac{{\bf1}_{X_{T_1}=2}}{\varphi_\omega^2}\bigg),
\eeqnn
where $V_i=\#\{k\in [0,T_1):X_k=i\}$. Using the shift invariance of $P$, we get
\beqnn
\int f \mbox{d}Q
&=&\sum_{i\le0}E_P\bigg(f(\theta^i\omega)\Big[\frac{1}{\varphi_\omega^1}E_\omega(V_i\cdot{\bf1}_{X_{T_1}=1})
                                             +\frac{1}{\varphi_\omega^2}E_\omega(V_i\cdot{\bf1}_{X_{T_1}=2})\Big]\bigg)\\
&=&\sum_{i\le0}E_P\bigg(f(\omega)\Big[\frac{1}{\varphi_{\theta^{-i}\omega}^1}E_{\theta^{-i}\omega}(V_i\cdot{\bf1}_{X_{T_1}=1})
                                     +\frac{1}{\varphi_{\theta^{-i}\omega}^2}E_{\theta^{-i}\omega}(V_i\cdot{\bf1}_{X_{T_1}=2})\Big]\bigg)\\
&=&E_P\bigg(f(\omega)\sum_{i\le0}\Big[\frac{1}{\varphi_{\theta^{-i}\omega}^1}E_{\theta^{-i}\omega}(V_i\cdot{\bf1}_{X_{T_1}=1})
                                     +\frac{1}{\varphi_{\theta^{-i}\omega}^2}E_{\theta^{-i}\omega}(V_i\cdot{\bf1}_{X_{T_1}=2})\Big]\bigg)\\
&=&E_P\bigg(f(\omega)\sum_{i\le0}\Big[E_{\theta^{-i}\omega}(V_i|X_{T_1}=1)+E_{\theta^{-i}\omega}(V_i|X_{T_1}=2)\Big]\bigg).
\eeqnn
Hence,
$$
\frac{\mbox{d}Q}{\mbox{d}P}=\sum_{i\le0}\Big[E_{\theta^{-i}\omega}(V_i|X_{T_1}=1)
+E_{\theta^{-i}\omega}(V_i|X_{T_1}=2)\Big].
$$
The calculation of $E_{\theta^{-i}\omega}(V_i|X_{T_1}=k),~(k=1,2)$ is based on the branching process.
Note that $V_i=A(i+1)+B(i+1)+C(i+1)+A(i)+B(i)=|U(i+1)|+A(i)+B(i)$. Then
\beqnn
E_\omega\Big(V_0\Big|X_{T_1}=1\Big)&=&E_\omega\Big[E_\omega\Big(|U(1)|+A(0)+B(0)~\Big|~U(1)\Big)\Big|X_{T_1}=1\Big]\\
&=&E_\omega\Big[|U(1)|+\frac{\alpha(0)}{1-\alpha(0)-\beta(0)}|U(1)|+\frac{\beta(0)}{1-\alpha(0)-\beta(0)}|U(1)|\Big|X_{T_1}=1\Big]\\
&=&\frac{1}{1-\alpha(0)-\beta(0)}\Big|E_\omega\Big(U(1)\Big|X_{T_1}=1\Big)\Big|\\
&=&\frac{1}{1-\alpha(0)-\beta(0)}\Big\langle(1,1,1), E_\omega\Big(U(1)\Big|X_{T_1}=1\Big)\Big\rangle.
\eeqnn
Similarly, for $i\le-1$,
\beqnn
E_\omega\Big(V_i\Big|X_{T_1}=1\Big)&=&E_\omega\Big[E_\omega\Big(V_i\Big|U(i+1)\Big)\Big|X_{T_1}=1\Big]\\
&=&\frac{1}{1-\alpha(i)-\beta(i)}\Big|E_\omega\Big(U(i+1)\Big|X_{T_1}=1\Big)\Big|\\
&=&\frac{1}{1-\alpha(i)-\beta(i)}\Big|E_\omega\Big(U(1)\Big|X_{T_1}=1\Big)\cdot N(0)\cdots N(i+1)\Big|\\
&=&\frac{1}{1-\alpha(i)-\beta(i)}\Big\langle(1,1,1),E_\omega\Big(U(1)\Big|X_{T_1}=1\Big)\cdot N(0)\cdots N(i+1)\Big\rangle.
\eeqnn
Hence, for $i\le-1$,
\beqnn
E_{\theta^{-i}\omega}\Big(V_i\Big|X_{T_1}=1\Big)=
\frac{1}{1-\alpha(0)-\beta(0)}\Big\langle (1,1,1), E_{\theta^{-i}\omega}\Big(U(1)\Big|X_{T_1}=1\Big)\cdot N(-i)\cdots N(1)\Big\rangle.
\eeqnn
By the same argument,
\beqnn
E_\omega\Big(V_0\Big|X_{T_1}=2\Big)=\frac{1}{1-\alpha(0)-\beta(0)}\Big\langle (1,1,1),E_\omega\Big(U(1)\Big|X_{T_1}=2\Big)\Big\rangle.
\eeqnn
And for $i\le-1$,
\beqnn
E_{\theta^{-i}\omega}\Big(V_i\Big|X_{T_1}=2\Big)=
\frac{1}{1-\alpha(0)-\beta(0)}\Big\langle(1,1,1),E_{\theta^{-i}\omega}\Big(U(1)\Big|X_{T_1}=2\Big)\cdot N(-i)\cdots N(1)\Big\rangle.
\eeqnn
By branching process, we have
\beqnn
& &E_\omega\Big(U(1)\Big|X_{T_1}=1\Big)\\
&=&(1,0,0)\cdot P_\omega(U(1)=(1,0,0)|X_{T_1}=1)\\
& &\qquad
+(0,1,0)\cdot P_\omega(U(1)=(0,1,0)|X_{T_1}=1)+(0,0,1)\cdot P_\omega(U(1)=(0,0,1)|X_{T_1}=1)\\
&=&(1,0,0)\cdot \frac{p_1(0)}{1-\alpha(0)-\beta(0)}\Big/\Big(\frac{p_1(0)+\gamma(0)}{1-\alpha(0)-\beta(0)}\Big)\\
& &\qquad\qquad\qquad\qquad+(0,1,0)\cdot \frac{\gamma(0)}{1-\alpha(0)-\beta(0)}\Big/\Big(\frac{p_1(0)+\gamma(0)}{1-\alpha(0)-\beta(0)}\Big)+
(0,0,1)\cdot 0\\
&=&\Big(\frac{p_1(0)}{p_1(0)+\gamma(0)},\frac{\gamma(0)}{p_1(0)+\gamma(0)},0\Big).
\eeqnn
And,
\beqnn
E_\omega\Big(U(1)\Big|X_{T_1}=2\Big)=(0,0,1).
\eeqnn
Consequently, we get that
\beqnn
\frac{\mbox{d} Q}{\text{d} P}&=&\sum_{i\le0}\bigg(\sum_{k=1,2}E_{\theta^{-i}\omega}\Big(V_i\Big|X_{T_1}=k\Big)\bigg)\\
&=&\frac{2+\bigg\langle(1,1,1),\sum_{i\ge1}\Big[\sum_{k=1,2}E_{\theta^{i}\omega}\Big(U(1)\Big|X_{T_1}=k\Big)\Big]\cdot N(i)\cdots N(1)\bigg\rangle}
{1-\alpha(0)-\beta(0)}\\
&=&\frac{2+\bigg\langle(1,1,1),\sum_{i\ge1}\Big(\frac{p_1(i)}{p_1(i)+\gamma(i)},\frac{\gamma(i)}{p_1(i)+\gamma(i)},1\Big)\cdot N(i)\cdots N(1)\bigg\rangle}
{1-\alpha(0)-\beta(0)}.
\eeqnn
\qed


\bigskip

\section{The law of large numbers\label{s4}}

\setcounter{equation}{0}

We are now ready to prove the law of large number with an {\it
explicit } drift based on the branching structure by the method of
``environment viewed from the particle".

\blemma\label{l4.1}
Under the law induced by $\overline{Q}\otimes
P_\omega$, the sequence $\{\overline{\omega}(n)\}$ is stationary and
ergodic.
\elemma

\noindent{\it Proof }. Since $\overline{Q}$ is an invariant measure
under the transition kernel $K$ by Theorem \ref{t1.2}, we obtain that the
process $\{\overline{\omega}(n)\}$ is stationary under
$\overline{Q}\otimes P_\omega$. For the ergodicity, the proof is
similar  as Zeitouni (\cite{Z04}, page 207, Corollary 2.1.25), we
omit the details. \qed

\medskip

\noindent{\it Proof of Theorem \ref{t1.3}}. The idea of the proof
for the LLN is the same as Zeitouni (\cite{Z04}, page 208, Theorem
2.1.9), however, we pay attention to the drift here.

Define the local drift at site $x$ in the environment $\omega$, as
$d(x,\omega)=E_\omega^x(X_1-x)$. The ergodicity of
$\{\overline\omega(i)\}$ under $\overline Q\otimes P_\omega$ implies
that: \beqnn \frac1n\sum_{k=0}^{n-1}d(X_k,\omega)
=\frac1n\sum_{k=0}^{n-1}d(0,\overline\omega(k))\longrightarrow \int
d(0,\omega)~\mbox{d}\overline Q, ~\mbox{as}~ n\to \infty,~\overline
Q\otimes P_\omega\mbox{-a.s.} \eeqnn On the other hand, \beqnn
X_n=\sum_{i=1}^n(X_i-X_{i-1})&=&\sum_{i=1}^n\Big(X_i-X_{i-1}-d(X_{i-1},
\omega)\Big)+\sum_{i=1}^nd(X_{i-1},\omega)\\
&:=&M_n+\sum_{i=1}^nd(X_{i-1},\omega).
\eeqnn
Under $P_\omega$, $M_n$ is a martingale, with $|M_n-M_{n-1}|\le3$ for $\omega\in\Omega$.
Hence, with $\mathcal{G}_n=\sigma(M_1,\cdots,M_n)$,
\beqnn
E_\omega(\mbox{e}^{\lambda M_n})&=&E_\omega\Big(\mbox{e}^{\lambda
M_{n-1}}E_\omega(\mbox{e}^{\lambda
(M_n-M_{n-1})}|\mathcal{G}_n)\Big)\\
&\le&E_\omega\Big(\mbox{e}^{\lambda
M_{n-1}}\mbox{e}^{3\lambda^2}\Big),
\eeqnn
and hence, iterating,
$E_\omega(\mbox{e}^{\lambda M_n})\le \mbox{e}^{3n\lambda^2}$ (this
is a version of Azuma's inequality, see \cite{DZ98}, Corollary
2.4.7). Then Chebyshev's inequality implies that
\beqnn
\frac{M_n}{n}\to 0,~\mathbb{P}\mbox{-a.s.}
\eeqnn
Hence,
\beqnn
\lim_{n\to\infty}\frac1nX_n=\int d(0,\omega)~\mbox{d}\overline Q=v_P,
\eeqnn
Observe that
\beqnn
v_P&=&\int d(0,\omega)\mbox{d}\overline Q=E_{\overline Q}(X_1)\\
&=&\frac{E_P\Big[\Pi(\omega)\Big(p_1(0)+2p_2(0)-p_{-1}(0)\Big)\Big]}{Q(\Omega)}\\
&=&\frac{E_P\bigg[\frac{p_1(0)+2p_2(0)-p_{-1}(0)}{1-\alpha(0)-\beta(0)}\Big(2+
\Big\langle(1,1,1),\sum_{i\ge1}\Big(\frac{p_1(i)}{p_1(i)+\gamma(i)},
\frac{\gamma(i)}{p_1(i)+\gamma(i)},1\Big)\cdot N(i)\cdots N(1)\Big)\Big\rangle\bigg]}
{\mathbb{E}(T_1|X_{T_1}=1)+\mathbb{E}(T_1|X_{T_1}=2)}.
\eeqnn
Moreover, by (\ref{2.2}) and the proof of the density of Theorem \ref{t1.2}, we have
\beqnn
\mathbb{E}(T_1|X_{T_1}=1)&=&E_P\bigg[E_\omega\Big(1+\Big\langle(2,2,1), \sum_{i\le0}U(i)\Big\rangle\Big|X_{T_1}=1\Big)\bigg]\\
&=&E_P\bigg[1+\Big\langle (2,2,1), \sum_{i\le0} E_\omega \Big(U(i)\Big|X_{T_1}=1\Big)\Big\rangle \bigg]\\
&=&E_P\bigg[1+\Big\langle (2,2,1), \sum_{i\le0} E_\omega \Big(U(1)\Big|X_{T_1}=1\Big)\cdot N(0)\cdots N(i)\Big\rangle \bigg]\\
&=&E_P\bigg[1+\Big\langle (2,2,1), \sum_{i\le0} \Big( \frac{p_1(0)}{p_1(0)+\gamma(0)},
\frac{\gamma(0)}{p_1(0)+\gamma(0)},0\Big)\cdot N(0)\cdots N(i)\Big\rangle \bigg],
\eeqnn
and by the same argument,
\beqnn
\mathbb{E}(T_1|X_{T_1}=2)=
E_P\bigg[1+\Big\langle (2,2,1), \sum_{i\le0} \Big(0,0,1\Big)\cdot N(0)\cdots N(i)\Big\rangle \bigg].
\eeqnn
Therefore,
\beqnn
& &\mathbb{E}(T_1|X_{T_1}=1)+\mathbb{E}(T_1|X_{T_1}=2)\\
&=&E_P\bigg[2+\Big\langle (2,2,1), \sum_{i\le0} \Big(\frac{p_1(0)}{p_1(0)+\gamma(0)},\frac{\gamma(0)}{p_1(0)+\gamma(0)},1\Big)
\cdot N(0)\cdots N(i)\Big\rangle \bigg]\\
&=&E_P\bigg[2+\Big\langle (2,2,1), \sum_{i\ge0} \Big(\frac{p_1(i)}{p_1(i)+\gamma(i)},
\frac{\gamma(i)}{p_1(i)+\gamma(i)},1\Big)\cdot N(i)\cdots N(0)\Big\rangle \bigg],
\eeqnn
where the last equation holds by the invariance of $P$. Consequently, the drift in this case can be written as
\beqnn
v_P=\frac{E_P\bigg[\frac{p_1(0)+2p_2(0)-p_{-1}(0)}{1-\alpha(0)-\beta(0)}
\Big(2+\sum_{i\ge1}\Big\langle(1,1,1),
\Big(\frac{p_1(i)}{p_1(i)+\gamma(i)},\frac{\gamma(i)}{p_1(i)+\gamma(i)},1\Big)\cdot N(i)\cdots N(1)\Big\rangle\Big)\bigg]}
{E_P\bigg[2+\Big\langle (2,2,1),
\sum_{i\ge0} \Big(\frac{p_1(i)}{p_1(i)+\gamma(i)},\frac{\gamma(i)}{p_1(i)+\gamma(i)},1\Big)\cdot N(i)\cdots N(0)\Big\rangle \bigg]}.
\eeqnn
\qed

\noindent {\bf Remark} Consider the random walk in non-random
environment defined in Proposition \ref{p1.1}. Then the drift $v_P$
reduces to $E(X_1)=p_1+2p_2-q$. In fact, when the environment is not
random, the drift in Theorem \ref{t1.3} can be written as
\beqnn
v_P=(p_1+2p_2-q)\frac{\frac{1}{1-\alpha-\beta}\Big(2+\Big\langle(1,1,1),
\Big(\frac{p_1}{p_1+\gamma},\frac{\gamma}{p_1+\gamma},1\Big)\cdot\sum_{n\ge1}N^n\Big\rangle\Big)}
{2+\Big\langle (2,2,1),
\Big(\frac{p_1}{p_1+\gamma},\frac{\gamma}{p_1+\gamma},1\Big)\cdot\sum_{n\ge1}N^n\Big\rangle}.
\eeqnn
By (\ref{2.22}), we obtain
\beqnn
v_P&=&(p_1+2p_2-q)\\
& &\times\frac{\frac{1}{1-\alpha-\beta}\Big[2+\frac{1}{(p_1+\gamma)(1-2\alpha-3\beta)}\Big\langle(1,1,1),
\Big(\alpha(2p_1+3\gamma),~\beta(2p_1+3\gamma),~2\beta p_1+(1-2\alpha)\gamma\Big)\Big\rangle\Big]}
{2+\frac{1}{(p_1+\gamma)(1-2\alpha-3\beta)}\Big\langle(2,2,1),
\Big(\alpha(2p_1+3\gamma),~\beta(2p_1+3\gamma),~2\beta p_1+(1-2\alpha)\gamma\Big)\Big\rangle}\\
&=&(p_1+2p_2-q)\times\frac{~~\frac{2p_1+3\gamma}{(p_1+\gamma)(1-2\alpha-3\beta)}~~}{~~\frac{2p_1+3\gamma}{(p_1+\gamma)(1-2\alpha-3\beta)}~~}\\
&=&p_1+2p_2-q,
\eeqnn
as it should be.
\qed


\bigskip

\section{The general bounded jump case: $(1,R)$-RWRE\label{s5}}

\setcounter{equation}{0}

In this section, we consider $(1,R)$-RWRE, in which the possible jumps to
the right are $1,2,\cdots, R$, where $R$ is some fixed positive
integer. In this case, the environment is an element
$$
\omega=\{(q(\omega)_z,p_1(\omega)_z,\cdots,p_R(\omega)_z):z\in\mathbb{Z}\}
\in\mathcal{M}(\widetilde\Lambda)^\mathbb{Z}=:\widetilde\Omega,
$$
where $\widetilde\Lambda=\{-1,1,\cdots,R\}$ is the set of all possible jumps.
Following the idea in section \ref{s2}, we can define a branching process
with $(1+2+\cdots+R)$-type species. In fact, for a jump from $i$ to $i-1$,
there are $(1+2+\cdots+R)$ crossing-back ways, i.e., jumping from $i-k_1$ to $i$ for
$k_1\in\{1,2,\cdots,R\}$, jumping from $i-k_2$ to $i+1$ for $k_2\in\{1,2,\cdots,R-1\}$, $\cdots$,
jumping from $i-k_R$ to $i+R-1$ for $k_R=1$.

  First, we need the following lemma dealing with exit
probabilities for this general case.

\blemma\label{l5.1} For $n\ge2$ and $1\le j\le R$, we have \beqnn
& &P_\omega^{i}[(-n,i),i+j]=\frac{\langle e_1,
[\widetilde{M}(i)+\cdots +\widetilde{M}(-n)\cdots
\widetilde{M}(i)](e_{j}-e_{j+1})\rangle} {1+\langle e_1,
[\widetilde{M}(i)+\cdots +\widetilde{M}(-n)\cdots
\widetilde{M}(i)]e_1\rangle}, \eeqnn with $e_{R+1}=0$ and \beqnn
\widetilde {M}(i):= \left( \begin{array}{cccc}
\frac{p_1(i)+\cdots+p_R(i)}{q(i)} &  \cdots & \frac{p_{R-1}(i)+p_R(i)}{q(i)} & \frac{p_R(i)}{q(i)}\\
1 & \cdots & 0 &0 \\
\vdots& \ddots & \vdots &\vdots\\
0&\cdots&1&0
\end{array} \right).
\eeqnn \elemma

The proof is the same as Lemma \ref{l2.1}. \qed

\

Define for $i\le0$,
$$
U(i)=\Big(U_1(i),\cdots,U_R(i),U_{R+1}(i),\cdots,U_{R+R-1}(i),\cdots,U_{1+2+\cdots+R}(i)\Big),
$$
where, for $k=1,2,\cdots,(1+2+\cdots+R)$, $U_k(i)$ is the number of steps from $i$ to $i-1$ before time $T_1$ with crossing-back from
$i-1$ to $i$, $\cdots$, $i-R$ to $i$; $i-1$ to $i+1$, $\cdots$, $i-R+1$ to $i+1$; $\cdots$;
$i-1$ to $i+R-1$; respectively. For probabilities, let
\beqnn
& &p_{(1)}(i)=P_\omega^i[\mbox{the walk jumping from $i$ to $i-1$, and finally jumping back } \\
& &\qquad\qquad\quad\mbox{ from $i-1$ to $i$ before $T_1$}],\\
& &\quad~\vdots\\
& &p_{(R)}(i)=P_\omega^i[\mbox{the walk jumping from $i$ to $i-1$, and finally jumping back } \\
& &\qquad\qquad\quad\mbox{ from $i-R$ to $i$ before $T_1$}],\\
& &p_{(R+1)}(i)=P_\omega^i[\mbox{the walk jumping from $i$ to $i-1$, and finally jumping back } \\
& &\qquad\qquad\quad\mbox{ from $i-1$ to $i+1$ before $T_1$}],\\
& &\quad~\vdots\\
& &p_{(R+R-1)}(i)=P_\omega^i[\mbox{the walk jumping from $i$ to $i-1$, and finally jumping back } \\
& &\qquad\qquad\quad\mbox{ from $i-R+1$ to $i+1$ before $T_1$}],\\
& &\quad~\vdots\\
& &p_{(1+2+\cdots+R)}(i)=P_\omega^i[\mbox{the walk jumping from $i$ to $i-1$, and finally jumping back } \\
& &\qquad\qquad\quad\mbox{ from $i-1$ to $i+R-1$ before $T_1$}].
\eeqnn
Expressions of these probabilities can be calculated by
using the exit probabilities as this is done in section \ref{s2}. Firstly,
we have
\beqnn
p_{(1)}+\cdots+p_{(R)}
&=&P_\omega^i[\mbox{the walk jumping from $i$ to $i-1$, and finally }\\
& &\qquad\mbox{ jumping back to $i$ before $T_1$}],\\
p_{(R+1)}+\cdots+p_{(R+R-1)}
&=&P_\omega^i[\mbox{the walk jumping from $i$ to $i-1$, and finally }\\
& &\qquad\mbox{ jumping back to $i+1$ before $T_1$}],\\
\vdots~~& &\\
p_{((1+2+\cdots+R)-2)}(i)+p_{((1+2+\cdots+R)-1)}(i)
&=&P_\omega^i[\mbox{the walk jumping from $i$ to $i-1$, and finally }\\
& &\qquad\mbox{jumping back to $i+R-2$ before $T_1$}],\\
p_{(1+2+\cdots+R)}(i)
&=&P_\omega^i[\mbox{the walk jumping from $i$ to $i-1$, and finally }\\
& &\qquad\mbox{jumping back to $i+R-1$ before $T_1$}].
\eeqnn
Recall the exit probabilities, we obtain
\beqlb\label{5.1}
p_{(1)}+\cdots+p_{(R)}&=&q(i)\cdot P_\omega^{i-1}[(-\infty,i-1),i],\nonumber\\
p_{(R+1)}+\cdots+p_{(R+R-1)}&=&q(i)\cdot P_\omega^{i-1}[(-\infty,i-1),i+1],\nonumber\\
\vdots~~& &\nonumber\\
p_{((1+2+\cdots+R)-2)}(i)+p_{((1+2+\cdots+R)-1)}(i)&=&q(i)\cdot P_\omega^{i-1}[(-\infty,i-1),i+R-2],\\
p_{(1+2+\cdots+R)}(i)&=&q(i)\cdot P_\omega^{i-1}[(-\infty,i-1),i+R-1].\nonumber
\eeqlb
Observe that
\beqnn
p_{((1+2+\cdots+R)-2)}(i):p_{((1+2+\cdots+R)-1)}(i)=p_{R-1}(i-1):p_{(1+2+\cdots+R)}(i-1).
\eeqnn
Thus by (\ref{5.1}), we get
\beqnn
& &p_{((1+2+\cdots+R)-2)}(i)=q(i)\cdot P_\omega^{i-1}[(-\infty,i-1),i+R-2]
\cdot\frac{p_{R-1}(i-1)}{p_{R-1}(i-1)+p_{(1+2+\cdots+R)}(i-1)},\\
& &p_{((1+2+\cdots+R)-1)}(i)=q(i)\cdot P_\omega^{i-1}[(-\infty,i-1),i+R-2]
\cdot\frac{p_{(1+2+\cdots+R)}(i-1)}{p_{R-1}(i-1)+p_{(1+2+\cdots+R)}(i-1)}.
\eeqnn
Using the same argument, we can get the expression of all the probabilities:
\beqnn
p_{(1+2+\cdots+R)}(i)&=&q(i)\cdot P_\omega^{i-1}[(-\infty,i-1),i+R-1],\\
p_{((1+2+\cdots+R)-2)}(i)&=&q(i)\cdot P_\omega^{i-1}[(-\infty,i-1),i+R-2]\cdot\frac{p_{R-1}(i-1)}{p_{R-1}(i-1)+p_{(1+2+\cdots+R)}(i-1)},\\
p_{((1+2+\cdots+R)-1)}(i)&=&q(i)\cdot P_\omega^{i-1}[(-\infty,i-1),i+R-2]\cdot\frac{p_{(1+2+\cdots+R)}(i-1)}{p_{R-1}(i-1)+p_{(1+2+\cdots+R)}(i-1)},\\
\vdots\quad& &\\
p_{(R+1)}(i)&=&q(i)\cdot P_\omega^{i-1}[(-\infty,i-1),i+1]\\
& &\qquad\qquad~\times\frac{p_2(i-1)}{p_2(i-1)+p_{(R+R-1+1)}(i-1)+\cdots+p_{(R+R-1+R-2)}(i-1)},\\
\vdots\quad& &\\
p_{(R+R-1)}(i)&=&q(i)\cdot P_\omega^{i-1}[(-\infty,i-1),i+1]\\
& &\qquad\qquad~\times\frac{p_{(R+R-1+R-2)}(i-1)}{p_2(i-1)+p_{(R+R-1+1)}(i-1)+\cdots+p_{(R+R-1+R-2)}(i-1)},\\
p_{(1)}(i)&=&q(i)\cdot P_\omega^{i-1}[(-\infty,i-1),i]\cdot\frac{p_1(i-1)}{p_1(i-1)+p_{(R+1)}(i-1)+\cdots+p_{(R+R-1)}(i-1)},\\
p_{(2)}(i)&=&q(i)\cdot P_\omega^{i-1}[(-\infty,i-1),i]\cdot\frac{p_{(R+1)}(i-1)}{p_1(i-1)+p_{(R+1)}(i-1)+\cdots+p_{(R+R-1)}(i-1)},\\
\vdots\quad& &\\
p_{(R)}(i)&=&q(i)\cdot P_\omega^{i-1}[(-\infty,i-1),i]\cdot\frac{p_{(R+R-1)}(i-1)}{p_1(i-1)+p_{(R+1)}(i-1)+\cdots+p_{(R+R-1)}(i-1)}.
\eeqnn

\medskip

Now we are ready to give the theorem about $U(i)$.

\btheorem\label{t5.1}
Assume $X_n\to \infty,~\mathbb{P}$-a.s.. Then $\Big(U(i)\Big)_{i\le0}$ is an inhomogeneous
multitype branching process in $\mathbb{R}^{1+2+\cdots+R}$ with immigration
\beqnn
& &U(1)=e_1,\quad \mbox{with probability }\quad \frac{p_1(0)}{1-p_{(1)}(0)-\cdots-p_{(R)}(0)},\\
& &U(1)=e_2,\quad \mbox{with probability }\quad \frac{p_{(R+1)}(0)}{1-p_{(1)}(0)-\cdots-p_{(R)}(0)},\\
& &\quad \vdots\\
& &U(1)=e_R,\quad \mbox{with probability }\quad \frac{p_{(R+R-1)}(0)}{1-p_{(1)}(0)-\cdots-p_{(R)}(0)},\\
& &U(1)=e_{R+1},\quad \mbox{with probability }\quad \frac{p_2(0)}{1-p_{(1)}(0)-\cdots-p_{(R)}(0)},\\
& &\quad \vdots\\
& &U(1)=e_{2+3+\cdots+R},\quad \mbox{with probability }\quad \frac{p_{(1+2+\cdots+R)}(0)}{1-p_{(1)}(0)-\cdots-p_{(R)}(0)},\\
& &U(1)=e_{1+2+\cdots+R},\quad \mbox{with probability }\quad
\frac{p_R(0)}{1-p_{(1)}(0)-\cdots-p_{(R)}(0)},
\eeqnn
and the following offspring distribution:
\beqnn
& &P_\omega\Big(U(0)=(u_1,\cdots,u_R,0,\cdots,0),U(1)=e_1\Big||U(1)|=1\Big)\\
& &\qquad\qquad\qquad=\frac{(u_1+\cdots+u_R)!}{u_1!\cdots u_R!}p_{(1)}^{u_1}(0)\cdots p_{(R)}^{u_R}(0)p_1(0),\\
& &P_\omega\Big(U(0)=(u_1,\cdots,u_R,1,0,\cdots,0),U(1)=e_2\Big||U(1)|=1\Big)\\
& &\qquad\qquad\qquad=\frac{(u_1+\cdots+u_R)!}{u_1!\cdots u_R!}p_{(1)}^{u_1}(0)\cdots p_{(R)}^{u_R}(0)p_{(R+1)}(0),\\
& &\quad\vdots\\
& &P_\omega\Big(U(0)=(u_1,\cdots,u_R,0\cdots,0,\stackrel{\text{(R+R-1)th}}{1},0,\cdots,0),U(1)=e_R\Big||U(1)|=1\Big)\\
& &\qquad\qquad\qquad=\frac{(u_1+\cdots+u_R)!}{u_1!\cdots u_R!}p_{(1)}^{u_1}(0)\cdots p_{(R)}^{u_R}(0)p_{(R+R-1)}(0),\\
& &P_\omega\Big(U(0)=(u_1,\cdots,u_R,0,\cdots,0),U(1)=e_{R+1}\Big||U(1)|=1\Big)\\
& &\qquad\qquad\qquad=\frac{(u_1+\cdots+u_R)!}{u_1!\cdots u_R!}p_{(1)}^{u_1}(0)\cdots p_{(R)}^{u_R}(0)p_2(0),\\
& &\quad\vdots\\
& &P_\omega\Big(U(0)=(u_1,\cdots,u_R,0,\cdots,0,1),U(1)=e_{2+3+\cdots+R}\Big||U(1)|=1\Big)\\
& &\qquad\qquad\qquad=\frac{(u_1+\cdots+u_R)!}{u_1!\cdots u_R!}p_{(1)}^{u_1}(0)\cdots p_{(R)}^{u_R}(0)p_{(1+2+\cdots+R)}(0),\\
& &P_\omega\Big(U(0)=(u_1,\cdots,u_R,0,\cdots,0),U(1)=e_{1+2+\cdots+R}\Big||U(1)|=1\Big)\\
& &\qquad\qquad\qquad=\frac{(u_1+\cdots+u_R)!}{u_1!\cdots
u_R!}p_{(1)}^{u_1}(0)\cdots p_{(R)}^{u_R}(0)p_R(0).
\eeqnn
Furthermore, for $i\le0$, the offspring mean matrix of the
$(-i+1)$-th generation is:
\beqnn
\widetilde{N}(i)=(\widetilde{N}_1(i)~\widetilde{N}_2),
\eeqnn
where
\beqnn
\widetilde{N}_1(i)=\left(\begin{array}{ccc}
\frac{p_{(1)}(i)}{1-p_{(1)}(i)-\cdots-p_{(R)}(i)}& \cdots & \frac{p_{(R)}(i)}{1-p_{(1)}(i)-\cdots-p_{(R)}(i)}\\
\vdots&\vdots&\vdots \\
\frac{p_{(1)}(i)}{1-p_{(1)}(i)-\cdots-p_{(R)}(i)}& \cdots & \frac{p_{(R)}(i)}{1-p_{(1)}(i)-\cdots-p_{(R)}(i)}
\end{array}\right)_{(1+2+\cdots+R)\times R},
\eeqnn
and
\beqnn
\widetilde{N}_2=\left(\begin{array}{llll}
Z_{1,R-1}   & Z_{1, R-2}   & \cdots   & Z_{1,1}\\
I_{R-1}     & Z_{R-1,R-2}  & \cdots   & Z_{R-1,1}\\
Z_{1,R-1}   & Z_{1,R-2}    & \cdots   & Z_{1,1}\\
Z_{R-2,R-1} & I_{R-2}      & \cdots   & Z_{R-2,1}\\
Z_{1,R-1}    & Z_{1,R-2}    & \cdots   & Z_{1,1}\\
~~~\vdots   & ~~~\vdots    & ~\vdots  & ~~\vdots\\
Z_{1,R-1}   & Z_{1,R-2}    & \cdots   & I_1\\
Z_{1,R-1}   & Z_{1,R-2}    & \cdots   & Z_{1,1}
\end{array}\right)_{(1+2+\cdots+R)\times(1+2+\cdots+R-1)},
\eeqnn
in which $Z_{m,n}$ is the zero matrix of dimension $m\times n$, and
$I_m$ is the identity matrix of dimension $m\times m$.
\etheorem

\medskip

\noindent{\bf Remark} Following the ideas in section \ref{s2},
Theorem \ref{t5.1} can be proved analogously. We omit the details.

\

\

\noindent{\bf Acknowledgements} The authors would like to thank
Professor Zenghu Li,  Hongyan Sun and  Huanming Wang for the
stimulating discussions.


\end{document}